\documentclass[11pt,a4paper]{article}
\usepackage[latin1]{inputenc}
\usepackage{amsmath}
\usepackage{amsfonts}
\usepackage{amssymb}
\usepackage{makeidx}
\numberwithin{equation}{section}
\usepackage[width=12.00cm, height=19.50cm, left=2.50cm, right=2.50cm, top=2.50cm, bottom=2.cm]{geometry}
\author{Kuldeep Singh Gehlot}
\title{Properties of Ultra Gamma Function}
\begin{document}
\maketitle
\begin{center}
Government College Jodhpur,\\
JNV University Jodhpur, Rajasthan, India-306401.\\
Email: drksgehlot@rediffmail.com
\end{center}
\section*{Abstract}
In this paper we study the integral of type
 \[_{\delta,a}\Gamma_{\rho,b}(x) =\Gamma(\delta,a;\rho,b)(x)=\int_{0}^{\infty}t^{x-1}e^{-\frac{t^{\delta}}{a}-\frac{t^{-\rho}}{b}}dt.\]
Different authors called this integral by different names like ultra gamma function, generalized gamma function, Kratzel integral, inverse Gaussian integral, reaction-rate probability integral, Bessel integral etc.
We prove several identities and recurrence relation of above said integral, we called this integral as Four Parameter Gamma Function. Also we evaluate relation between Four Parameter Gamma Function, p-k Gamma Function and Classical Gamma      Function. With some conditions we can evaluate Four Parameter Gamma Function in term of Hypergeometric function.\\\\
\textbf{Mathematics Subject Classification :} 33B15.\\\\
\textbf{Keywords:} Four Parameter Gamma Function, Ultra Gamma Function, Two Parameter Gamma Function, Two Parameter Pochhammer Symbol. 
\section{Introduction}
The main aim of this paper is to introduce Four Parameter Gamma Function in the form,
\begin{equation}
_{\delta,a}\Gamma_{\rho,b}(x) =\Gamma(\delta,a;\rho,b)(x)=\int_{0}^{\infty}t^{x-1}e^{-\frac{t^{\delta}}{a}-\frac{t^{-\rho}}{b}}dt.
\end{equation}
where $ x\in C/\delta Z^{-} ; \delta, \rho, a,b \in R^{+}-\lbrace 0 \rbrace $ and  $ Re(x-\rho n)> 0, n\in N.$ \\
Four Parameter Gamma Function is the deformation of the two parameter Gamma Function defined by [4], such that          
$_{\delta,a}\Gamma_{\rho,b}(x)\Rightarrow\: _{k,p}\Gamma_{0,b}(x)=e^{-\frac{1}{b}}\: _{p}\Gamma_{k}(x)$, as $\delta=k,a=p,\rho=0. $\\
And this Four Parameter Gamma Function is the deformation of the k-Gamma Function defined by [1], such that $_{\delta,a}\Gamma_{\rho,b}(x)\Rightarrow\: _{k,k}\Gamma_{0,b}(x)=e^{-\frac{1}{b}}\:\Gamma_{k}(x)$, as $\delta=a=k,\rho=0. $\\ 
 Also the Four Parameter Gamma Function is the deformation of the classical Gamma Function, such that
 $_{\delta,a}\Gamma_{\rho,b}(x)\Rightarrow\: _{1,1}\Gamma_{0,b}(x)=e^{-\frac{1}{b}}\: \Gamma(x)$, as $\delta =a=1,\rho=0. $\\
Throughout this paper Let $ C,R^{+}, Re(),Z^{-}$ and $ N $ be the sets of complex numbers, positive real numbers, real part of complex number, negative integer and natural numbers respectively. We use the notation and terminology of [2] and [3].\\ 
The p - k Gamma Function (i.e. Two Parameter Gamma Function), $_{p}\Gamma_{k}(x)$ is given by [4],
For $ x\in C/kZ^{-};  k,p \in R^{+}-\lbrace 0 \rbrace $ and $ Re(x)>0, n\in N, $ is 
\begin{equation}
_{p}\Gamma_{k}(x)=\frac{1}{k}\lim_{n\rightarrow \infty} \dfrac{n!p^{n+1}(np)^{\frac{x}{k}}}{_{p}(x)_{n+1,k}}.
\end{equation}  
 or
\begin{equation}
 _{p}\Gamma_{k}(x)=\frac{1}{k}\lim_{n\rightarrow \infty} \dfrac{n!p^{n+1}(np)^{\frac{x}{k}-1}}{_{p}(x)_{n,k}}.
\end{equation} 
And the integral representation of p - k Gamma Function is given by
\begin{equation}
_{p}\Gamma_{k}(x)=\int^{\infty}_{0}e^{-\frac{t^{k}}{p}}t^{x-1}dt.
 \end{equation}  
\section{Recurrence formulas and infinite products of Four Parameter Gamma Function, $ _{\delta,a}\Gamma_{\rho,b}(x)$ or $  \Gamma(\delta,a;\rho,b)(x).$}   

\textbf{Theorem 2.1} The relation between Four Parameter Gamma Function, p - k Gamma Function and Classical Gamma Function is given by 
\begin{equation}
\Gamma(\delta,a;\rho,b)(x)=\sum_{n=0}^{\infty}\dfrac{(-1)^{n}}{n!b^{n}}\:_{a}\Gamma_{\delta}(x-\rho n),
\end{equation}
\begin{equation}
\Gamma(\delta,a;\rho,b)(x)=\sum_{n=0}^{\infty}\dfrac{(-1)^{n}}{n!b^{n}}\dfrac{1}{\delta}\:\lim_{m\rightarrow \infty} \dfrac{m!a^{m+1}(ma)^{\frac{x-\rho n}{\delta}-1}}{_{a}(x-\rho n)_{m,\delta}},
\end{equation}
And 
\begin{equation}
\Gamma(\delta,a;\rho,b)(x)=\sum_{n=0}^{\infty}\dfrac{(-1)^{n}}{n!b^{n}}\dfrac{a^{(\frac{x-\rho n}{\delta})}}{\delta}\:\Gamma(\frac{x-\rho n}{\delta}).
\end{equation}
where $ x\in C/\delta Z^{-} ; \delta, \rho, a,b \in R^{+}-\lbrace 0 \rbrace $ and  $ Re(x-\rho n)> 0, n\in N.$ \\
Proof: Using the definition (1.1), we have
\[\Gamma(\delta,a;\rho,b)(x)=\int_{0}^{\infty}t^{x-1}e^{-\frac{t^{\delta}}{a}-\frac{t^{-\rho}}{b}}dt, \]
\[\Gamma(\delta,a;\rho,b)(x)=\sum_{n=0}^{\infty}\dfrac{(-1)^{n}}{n!b^{n}}\int_{0}^{\infty}t^{x-\rho n-1}e^{-\frac{t^{\delta}}{a}}dt,\]
using [4], equation (2.14), we have
\[\Gamma(\delta,a;\rho,b)(x)=\sum_{n=0}^{\infty}\dfrac{(-1)^{n}}{n!b^{n}}\:_{a}\Gamma_{\delta}(x-\rho n),\]
And by using [4], theorem (2.9), we have
\[\Gamma(\delta,a;\rho,b)(x)=\sum_{n=0}^{\infty}\dfrac{(-1)^{n}}{n!b^{n}}\dfrac{a^{(\frac{x-\rho n}{\delta})}}{\delta}\:\Gamma(\frac{x-\rho n}{\delta}).\]
By using equation (1.3) in (2.1), we get the result (2.2).\\
This completes the proof.\\\\
\textbf{Theorem 2.2} Given $ x\in C/\delta Z^{-} ; \delta, \rho,k,p, a,b \in R^{+}-\lbrace 0 \rbrace $ and  $ Re(x-\rho n)> 0, n\in N,$ then the following recurrence relations exists,
\begin{equation}
\Gamma(\delta,a;\rho,b)(x)= \dfrac{k}{\delta}\;\Gamma(k,a;\frac{k\rho}{\delta},b)(\frac{kx}{\delta}),
\end{equation}
\begin{equation}
\Gamma(\delta,a;\rho,b)(x)= \dfrac{k}{\delta}(\frac{a}{p})^{\frac{x}{\delta}}\;\Gamma(k,p;\frac{k\rho}{\delta},b(\frac{a}{p})^{\frac{\rho}{\delta}})(\frac{kx}{\delta}),
\end{equation}
\begin{equation}
\Gamma(\delta,a;\rho,b)(x)= (\dfrac{a}{p})^{\frac{x}{\delta}}\;\Gamma(\delta,p;\rho,b(\frac{a}{p})^{\frac{\rho}{\delta}})(x),
\end{equation}
 Proof: From equation (2.1), we have
 \[\Gamma(\delta,a;\rho,b)(x)=\sum_{n=0}^{\infty}\dfrac{(-1)^{n}}{n!b^{n}}\:_{a}\Gamma_{\delta}(x-\rho n),\] 
using result (2.11) of [4],we have
\[\Gamma(\delta,a;\rho,b)(x)=\sum_{n=0}^{\infty}\dfrac{(-1)^{n}}{n!b^{n}}\dfrac{k}{\delta}\:_{a}\Gamma_{k}(\frac{k(x-\rho n)}{\delta}),\]
using equation (2.14) of [4], we have
\[\Gamma(\delta,a;\rho,b)(x)=\dfrac{k}{\delta}\int_{0}^{\infty}t^{\frac{kx}{\delta}-1}e^{-\frac{t^{k}}{a}-\frac{t^{-\frac{k\rho}{\delta}}}{b}}dt,\]
\[\Gamma(\delta,a;\rho,b)(x)=\dfrac{k}{\delta}\Gamma(k,a;\frac{k\rho}{\delta},b)(\frac{kx}{\delta}).\]
Which completes the proof of (2.4).\\
Similarly we can prove the result (2.5) and (2.6). \\\\
\textbf{Theorem 2.3} Given $ x\in C/\delta Z^{-} ; \delta, \rho,a,b \in R^{+}-\lbrace 0 \rbrace $ and  $ Re(x-\rho n)> 0, n\in N,$ then we can represent the Four Parameter Gamma Function in the form of series,
\begin{equation}
\Gamma(\delta,a;\rho,b)(x)= \dfrac{1}{\delta}(a)^{\frac{x}{\delta}}\sum_{n=0}^{\infty}\dfrac{(-1)^{n}}{n!}(\dfrac{a^{-\frac{\rho}{\delta}}}{b})^{n}\prod_{m=1}^{\infty}\lbrace(1-\frac{1}{m})^{(\frac{x-n\rho}{\delta})}(1+\frac{x-n\rho}{m\delta})^{-1}\rbrace
\end{equation}
Proof: From equation (2.1) we have,
 \[\Gamma(\delta,a;\rho,b)(x)=\sum_{n=0}^{\infty}\dfrac{(-1)^{n}}{n!b^{n}}\:_{a}\Gamma_{\delta}(x-\rho n),\] 
 using equation (2.15) of [4], we have
 \[\Gamma(\delta,a;\rho,b)(x)= \dfrac{1}{\delta}(a)^{\frac{x}{\delta}}\sum_{n=0}^{\infty}\dfrac{(-1)^{n}}{n!}(\dfrac{a^{-\frac{\rho}{\delta}}}{b})^{n}\prod_{m=1}^{\infty}\lbrace(1-\frac{1}{m})^{(\frac{x-n\rho}{\delta})}(1+\frac{x-n\rho}{m\delta})^{-1}\rbrace.\]
 Which completes the proof.\\\\
\textbf{Theorem 2.4} Given $ x\in C/\delta Z^{-} ; \delta, \rho,a,b \in R^{+}-\lbrace 0 \rbrace $ and  $ Re(x-\rho n)> 0, n\in N,$ then the fundamental equation satisfied by Four Parameter Gamma Function is,
\begin{equation}
x\:\Gamma(\delta,a;\rho,b)(x)= \dfrac{\delta}{a}\:\Gamma(\delta,a;\rho,b)(x+\delta)-\dfrac{\rho}{b}\:\Gamma(\delta,a;\rho,b)(x-\rho)
\end{equation}
Proof: From equation (2.1), we have
\[\Gamma(\delta,a;\rho,b)(x)=\sum_{n=0}^{\infty}\dfrac{(-1)^{n}}{n!b^{n}}\:_{a}\Gamma_{\delta}(x-\rho n),\]
replace $x $ by $x+\delta$ we have,
\[\Gamma(\delta,a;\rho,b)(x+\delta)=\sum_{n=0}^{\infty}\dfrac{(-1)^{n}}{n!b^{n}}\:_{a}\Gamma_{\delta}(x-\rho n+\delta),\]
using [4], equation (2.23), we have
\[\Gamma(\delta,a;\rho,b)(x+\delta)=\sum_{n=0}^{\infty}\dfrac{(-1)^{n}}{n!b^{n}}\dfrac{(x-\rho n)a}{\delta}\:_{a}\Gamma_{\delta}(x-\rho n),\]
\[\Gamma(\delta,a;\rho,b)(x+\delta)= \dfrac{xa}{\delta}\:\Gamma(\delta,a;\rho,b)(x)-\sum_{n=1}^{\infty}\dfrac{(-1)^{n}a\rho}{(n-1)!b^{n}\delta}\:_{a}\Gamma_{\delta}(x-\rho n),\]
$  n $ can be replace by $ n+1 $, we have
\[\Gamma(\delta,a;\rho,b)(x+\delta)= \dfrac{xa}{\delta}\:\Gamma(\delta,a;\rho,b)(x)+ \dfrac{a\rho}{b\delta}\sum_{n=0}^{\infty}\dfrac{(-1)^{n}}{(n)!b^{n}}\:_{a}\Gamma_{\delta}(x-\rho n-\rho),\]
using equation (2.1),
\[\Gamma(\delta,a;\rho,b)(x+\delta)= \dfrac{xa}{\delta}\:\Gamma(\delta,a;\rho,b)(x)+ \dfrac{a\rho}{b\delta}\:\Gamma(\delta,a;\rho,b)(x-\rho).\]
Which completes the proof.
\section{ Hypergeometric Function representation of Four Parameter Gamma Function, $ _{\delta,a}\Gamma_{\rho,b}(x)$ or $  \Gamma(\delta,a;\rho,b)(x).$}
\textbf{Theorem 3.1} Given $ x\in C/\delta Z^{-} ; \delta, \rho,a,b \in R^{+}-\lbrace 0 \rbrace, Re(x-\rho n)> 0, n\in N$ and $\dfrac{\rho}{\delta}\in N, $ then we have,
\begin{equation}
\Gamma(\delta,a;-\rho,b)(x)= \dfrac{a^{\frac{x}{\delta}}\Gamma(\frac{x}{\delta})}{\delta}\:_{\frac{\rho}{\delta}}F_{0}[(\frac{x-r\delta-\delta}{\rho})_{r=1,2,...,\frac{\rho}{\delta}};-;-\dfrac{1}{b}(\dfrac{a\rho}{\delta})^{\frac{\rho}{\delta}}].
\end{equation}
Proof: Consider the equation (2.3),
\[\Gamma(\delta,a;-\rho,b)(x)=\sum_{n=0}^{\infty}\dfrac{(-1)^{n}}{n!b^{n}}\dfrac{a^{(\frac{x+\rho n}{\delta})}}{\delta}\:\Gamma(\frac{x+\rho n}{\delta}),\] 
\[\Gamma(\delta,a;-\rho,b)(x)=\sum_{n=0}^{\infty}\dfrac{(-1)^{n}}{n!b^{n}}\dfrac{a^{(\frac{x+\rho n}{\delta})}}{\delta}\:\dfrac{\Gamma(\frac{x+\rho n}{\delta})\Gamma(\frac{x}{\delta})}{\Gamma(\frac{x}{\delta})},\]
\[\Gamma(\delta,a;-\rho,b)(x)=\Gamma(\frac{x}{\delta})\:\sum_{n=0}^{\infty}\dfrac{(-1)^{n}}{n!b^{n}}\dfrac{a^{(\frac{x+\rho n}{\delta})}}{\delta}\:(\frac{x}{\delta})_{\frac{\rho}{\delta}n},\]
we know the generalized pochammer symbol is given,
\[(\alpha)_{rn}=r^{rn}\prod_{n=1}^{r}(\frac{\alpha+n-1}{r})_{n},\] 
then we have,
\[\Gamma(\delta,a;-\rho,b)(x)=\Gamma(\frac{x}{\delta})\dfrac{a^{\frac{x}{\delta}}}{\delta}\:\sum_{n=0}^{\infty}\dfrac{1}{n!}[-\frac{1}{b}(\frac{a\rho}{\delta})^{\frac{\rho}{\delta}}]^{n}\prod_{r=1}^{\dfrac{\rho}{\delta}}(\frac{x-r\delta-\delta}{\rho})_{n},\]
this will give the desired result.\\\\
\section*{References}
\label{1}[1]  Diaz, R. and Pariguan, E. On hypergeometric functions and Pochhammer k-symbol. Divulgaciones Mathematicas, Vol. 15 No. 2 (2007) 179-192.\\
\label{2}[2]  Earl D. Rainville, Special Function, The Macmillan Company, New york,1963.\\
\label{3}[3]  Erdelyi, A., Higher Transcendental  Function Vol. 1, McGraw-Hill Book Company, New York, 1953.\\
\label{4}[4]  Kuldeep Singh Gehlot, Two Parameter Gamma Function and it's Properties, arXiv:1701.01052v1 [math.CA]
3 Jan 2017.\\ 
\label{5}[5]  Mathai, A.M. Computable Representation of Ultra Gamma Integral, J. of Ramanujan Society of Math. and Math. 
Vol.5, No.2 (2016), pp. 01-08. 
\end{document}